\documentclass[12pt,oneside,english]{amsart}
\usepackage{amssymb,amsmath,latexsym,amsthm}

\usepackage[a4paper,top=3cm,bottom=3cm,left=3cm,right=3cm,marginparwidth=1.75cm]{geometry}
\usepackage{graphicx}

\usepackage{xcolor}

\newtheorem{theorem}{Theorem}[section]
\newtheorem{corollary}[theorem]{Corollary}
\newtheorem{lemma}[theorem]{Lemma}
\newtheorem{prop}[theorem]{Proposition}

\theoremstyle{definition}
\newtheorem{definition}[theorem]{Definition}

\newcommand{\PSH}{{\rm PSH}}
\newcommand{\MPSH}{{\rm MPSH}}

 \usepackage{hyperref}
\hypersetup{  
    pdftitle={Cegrell classes},    
    pdfauthor={Mohammed Salouf},     
    colorlinks=true,       
   linkcolor=black,          
    citecolor=black,        
    filecolor=black,      
    urlcolor=black}        



\author{Omar ALEHYANE}
\author{Fatima Zahra ASSILA}
\author{Mohammed SALOUF}
\date{\today}
\address{
Laboratoire d'Informatique, Mathématiques et leurs Applications (LIMA), 
Faculté des Sciences El Jadida, 
Université Chouaib Doukkali,
24000 El Jadida, Maroc.} 
\email{alehyane.o@ucd.ac.ma}
\email{assila.f@ucd.ac.ma}
\email{salouf.m@ucd.ac.ma}
\begin{document}

\title[Degenerate complex Monge-Amp\`ere equations]{The Range of the Monge-Ampère operator $(\omega + dd^c .)^n$ in bounded domains}

\begin{abstract}
    Let $\Omega$ be a bounded strictly pseudoconvex domain of $\mathbb{C}^n$. We solve degenerate complex Monge-Ampère equations of the form $(\omega + dd^c \varphi)^n = \mu$ in the generalized Cegrell classes $\mathcal{K}(\Omega,\omega,H)$, where $H \in \mathcal{E}(\Omega)$ is maximal, $\omega$ is a smooth real $(1,1)$-form defined in a neighborhood of $\bar\Omega$ and $\mu$ is a positive Radon measure. This generalizes the previous work of the last author \cite{Sal25} to the case of non-continuous functions $H$ and also to the case of measures $\mu$ which do not vanish on pluripolar sets. 
\end{abstract}
\maketitle
\section{Introduction}
Let $\Omega$ be a bounded strictly pseudoconvex domain of $\mathbb{C}^n$, and let $\omega$ be a smooth real $(1,1)$-form defined in a neighborhood of $\bar\Omega$. 
Given a positive Radon measure $\mu$ on $\Omega$, we aim to study the complex Monge-Ampère equation 
 \begin{equation}\label{Main eq}
     (\omega + dd^c \varphi)^n = \mu.
 \end{equation}
 The left hand side of the equation refers to the Monge-Ampère measure of $\varphi$,  it is well defined on $\mathcal{C}^2$-smooth functions by the exterior product of $(1,1)$-forms with continuous coefficients. 

 When $\omega = 0$, Bedford-Taylor \cite{BT76} defined the Monge-Ampère measure $(dd^c \varphi)^n$ as a positive Radon measure for every function $\varphi$ which is locally bounded and plurisubharmonic (we write psh for short). Then, they proved that \eqref{Main eq}  admits continuous solutions in the case $\mu = f dV$, $f \in L^p(dV)$ and $p>1$ \cite[Theorem D]{BT76}. Based on this result, Ko\l odziej \cite[Page 1322]{Kol95} established a subsolution theorem for \eqref{Main eq}:  he proved that the equation \eqref{Main eq} admits bounded solutions if $\mu$ is dominated by the Monge-Ampère measure of a bounded psh function (we refer to \cite[Section 4]{Kol05} for an alternative proof).
Later on, Cegrell defined a general class $\mathcal{E}(\Omega)$ of unbounded psh functions for which the Monge-Ampère operator $(dd^c.)^n$ is well defined and continuous with respect to decreasing sequences \cite[Definition 4.1]{Ceg04}. Then, in a series of papers \cite{Ceg98,Ceg04,Ceg08}, Cegrell introduced many subsets of $\mathcal{E}(\Omega)$ in which he characterized the range of the operator $(dd^c.)^n$. In \cite[Theorem 4.13]{ACCH09}, \AA hag, Cegrell, Czy{\.z} and Hiep established a general subsolution theorem for \eqref{Main eq}. Namely, they proved that the equation \eqref{Main eq} has a solution in a Cegrell class $\mathcal{K}(\Omega)$ if $\mu$ is dominated by the Monge-Ampère measure of a function belonging to $\mathcal{K}(\Omega)$, which was an interesting generalization of \cite[Theorem 6.2]{Ceg04}. An alternative proof of this result can be found in \cite{QH19}.

Assume now $\omega \neq 0$. Building on the work of Bedford-Taylor \cite{BT76}, Ko\l odziej-Nguyen defined the operator $(\omega + dd^c .)^n$ on locally bounded $\omega$-psh functions \cite[Pages 143-144]{KN15Phong}. They also proved that \eqref{Main eq} has a bounded solution in the case $\mu = f dV$, $f \in L^p(dV)$ and $p>1$ \cite[Corollary 0.4]{KN15Phong}. Furthermore, by assuming $\mu \leq (dd^c u)^n$ is  dominated by the Monge-Ampère measure of a bounded psh function, Ko\l odziej and Nguyen proved that \eqref{Main eq} has bounded solutions
  \cite[Theorem 3.1]{KN23a}  and \cite[Theorem 3.2]{KN23b}.  In a previous work \cite{Sal25}, the last author defined a generalization of the Cegrell classes as follows: if $\mathcal{K}(\Omega)$ is a Cegrell class and $H \in \mathcal{E}(\Omega)$ is maximal, then $\mathcal{K}(\Omega,\omega,H)$ is given by  
  $$ u \in  \mathcal{K}(\Omega,\omega, H) \Leftrightarrow u \in \PSH(\Omega,\omega) \; \text{and} \;  u + \rho \in \mathcal{K}(\Omega,H), $$
  where $\rho$ is a $\mathcal{C}^2$-smooth function in a neighborhood of $\bar \Omega$ which satisfies $\omega \leq dd^c \rho$ (see Definition \ref{Def Gen Ceg class}). The last author proved that the operator $(\omega + dd^c .)^n$ is well defined on the generalized Cegrell classes $\mathcal{K}(\Omega,\omega,H)$ and that it is continuous with respect to decreasing sequences \cite[Theorem 4.10]{Sal25}. Moreover, he established solutions to \eqref{Main eq} in $\mathcal{K}(\Omega,\omega,H)$ under the conditions that $H$ is continuous on $\bar\Omega$ and $\mu$ vanishes on pluripolar sets and is dominated by the Monge-Ampère measure of a function in $\mathcal{K}(\Omega)$ \cite[Theorem 5.8]{Sal25}. 

  Here we prove the following. 
   \begin{theorem}[Main theorem]
  Let $H \in \mathcal{E}(\Omega)$ be a maximal function.  If $\mu \leq (dd^c u)^n$ is dominated by the Monge-Ampère measure of a function $u \in \mathcal{K}(\Omega)$ then \eqref{Main eq} admits a solution $\varphi \in \mathcal{K}(\Omega,\omega,H)$. 
\end{theorem}
We highlight that, in general, solutions to \eqref{Main eq} are not uniquely determined as demonstrated in \cite[Example 3.4]{Zer97}. We refer the reader to Theorem \ref{thm : uniqueness} for a partial result which is a generalization of \cite[Theorem 3.5]{ACLR25}. When the measure $\mu$ vanishes on pluripolar sets,  we know by \cite[Corollary 5.7]{Sal25} that solutions in the last theorem are uniquely determined.

To prove the last theorem, we first write $\mu$ as
$$ \mu = f(dd^c\psi)^n + (dd^cw)^n,   $$
where $w \in \mathcal{K}(\Omega,H)$, $\psi$ is in the Cegrell class $\mathcal{E}_0(\Omega)$ and $f \in L^1((dd^c \psi)^n)$. The existence of this decomposition is justified by \cite[Theorem 5.11]{Ceg04} and \cite[Theorem 4.13]{ACCH09}. Then, motivated by the proof of \cite[Theorem 4.13]{ACCH09}, we consider a generalized quasi-psh envelopes of the form 
\begin{equation}\label{quasi psh env}
    \varphi := \sup \{h \in \mathcal{E}(\Omega,\omega) : (\omega + dd^c h)^n \geq f (dd^c \psi)^n \; \text{and} \; \varphi + \rho \leq w \},
\end{equation}
where $\rho \in \mathcal{C}^2(\bar\Omega)$ is such that $\rho = 0$ on $\partial \Omega$ and $\omega \leq dd^c \rho$. The objective now is to prove that $\varphi$ is the solution in the last theorem. A major difficulty here is that, unlike the case $\omega=0$, the total masses of the Monge-Ampère measures are not monotone (see \cite[Example 5.1]{Sal25}). To overcome this difficulty, we studied quasi-psh envelopes of  the form \eqref{quasi psh env} in the spirit of \cite[Section 9]{BT82}. This allows us to prove that $\varphi \in \mathcal{K}(\Omega,\omega,H)$ and that the Monge-Ampère measure of $\varphi$ is equal to 
$f(dd^c \psi)^n$ on $\{\varphi>-\infty\}$. Then, by observing that 
$$ (\omega + dd^c \varphi)^n = (dd^c (\varphi+\rho))^n\; \; \text{on} \; \{\varphi = -\infty\}, $$
we were able to prove that the Monge-Ampère measure of $\varphi$ is equal to 
$(dd^c w)^n$ on $\{\varphi = -\infty\}$, hence $\varphi$ is the solution in the last theorem.

The article is organized as follows. In Section \ref{sec 2}, we recall some basic facts about quasi-psh functions and Cegrell classes, and we study a generalized quasi-psh envelopes. The main theorem is proved in Section \ref{sec: solving MAE}. 
\section*{Acknowledgment}
The authors thank Chinh H. Lu for helpful comments on an earlier draft and for suggesting Theorem \ref{thm : uniqueness}. 
\section{Generalized Cegrell classes and Quasi-psh envelopes}\label{sec 2}
Throughout this paper, we will denote by $\Omega$ a bounded strictly pseudoconvex domain of $\mathbb{C}^n$. We let $\omega$ be a smooth real $(1,1)$-form defined in a neighborhood of $\bar\Omega$, and we fix $\rho \in \PSH(\Omega) \cap \mathcal{C}^\infty(\bar\Omega)$ such that $\rho = 0$ on $\partial \Omega$ and $\omega \leq dd^c \rho$.

In this section, we first recall some basic facts about quasi-psh functions and Cegrell classes. Then, we study generalized quasi-psh envelopes. As we will see in Section \ref{sec: solving MAE}, these envelopes will be useful in solving degenerate complex Monge-Ampère equations. 
\subsection{Cegrell classes}
Let us denote by $\PSH(\Omega)$ the set of plurisubharmonic functions on $\Omega$. In \cite{Ceg98,Ceg04}, Cegrell defined the following subsets of $\PSH(\Omega)$: 
$$ \mathcal{E}_0(\Omega) = \{ u \in \PSH(\Omega)\cap L^{\infty}(\Omega) : u = 0 \; \text{on} \; \partial{\Omega} \; \text{and} \; \int_\Omega (dd^c u)^n < +\infty    \},
$$
$$ \mathcal{F}(\Omega) = \{ u \in \PSH(\Omega) : \exists u_j \in \mathcal{E}_0(\Omega), \; u_j \searrow u \; \text{and} \; \sup_j \int_\Omega (dd^c u_j)^n < +\infty    \}, $$ 
$$ \mathcal{E}(\Omega) = \{ u \in \PSH(\Omega) : \forall D \Subset \Omega, \exists u_D \in \mathcal{F}(\Omega) \; \text{such that} \; u_D = u \; \text{on} \; D \}, $$
and for $p>0$,
$$ \mathcal{E}_p(\Omega) = \{ u \in \PSH(\Omega) : \exists (u_j)_j \subset \mathcal{E}_0(\Omega), \; u_j \searrow u \; \; \text{and} \; \sup_j \int_\Omega   |u_j|^p (dd^c u_j)^n <+\infty  \}. $$
In \cite{Ceg04}, Cegrell proved that $\mathcal{E}(\Omega)$ contains all the psh functions that are  locally bounded and negative on $\Omega$. He also proved that the Monge-Ampère operator $(dd^c .)^n$ is well defined on $\mathcal{E}(\Omega)$ and continuous with respect to decreasing sequences, generalizing previous results of Bedford-Taylor \cite{BT76} which treats the case of locally bounded functions. Moreover, \cite[Theorem 4.5]{Ceg04} shows that $\mathcal{E}(\Omega)$ is the largest set of negative psh functions for which these two properties hold. 
 
Recall that a function $H \in \mathcal{E}(\Omega)$ is maximal if and only if it satisfies $(dd^cH)^n=0$, we let $\MPSH(\Omega)$ denote the set of maximal functions in $\Omega$. In \cite{Ceg08}, Cegrell defined the class $\mathcal{N}(\Omega)$ by the set of functions $u \in \mathcal{E}(\Omega)$ for which the smallest maximal function 
majorant of $u$ is identical to zero.  In other words, 
$$ \mathcal{N}(\Omega) = \{ u \in \mathcal{E}(\Omega) : \forall H \in \mathcal{E}(\Omega)\cap \MPSH(\Omega), \; u \leq H \leq 0 \Rightarrow H = 0 \}.   $$
We have 
$$ \mathcal{E}_0(\Omega) \subset \mathcal{F}(\Omega) \cap \mathcal{E}_p(\Omega) \subset \mathcal{F}(\Omega) \cup \mathcal{E}_p(\Omega)  \subset \mathcal{N}(\Omega) \subset \mathcal{E}(\Omega). $$
In particular, the Monge-Ampère operator $(dd^c.)^n$ is well defined on all these classes.

In \cite{Sal25}, the last author extended the results of Cegrell to the study of the operator $(\omega + dd^c .)^n$.
Recall that a function $u \in L^1_{loc}(\Omega)$ is quasi-psh on $\Omega$ if we can write $u = u_1 + u_2$ with $u_1 \in \mathcal{C}^\infty(\Omega)$ and $u_2 \in \PSH(\Omega)$. We denote by $\PSH(\Omega,\omega)$ the set of quasi-psh functions $u$ which satisfy $\omega + dd^c u \geq 0$ in the weak sense of currents; functions in  
$\PSH(\Omega,\omega)$ are called $\omega$-psh functions. It is easy to see that $\PSH(\Omega,\omega) = \PSH(\Omega)$ when $\omega=0$ and that $\PSH(\Omega) \subset \PSH(\Omega,\omega)$ when $\omega \geq 0$. Also, by our choice of the function $\rho$, we have 
$$ u + \rho \in \PSH(\Omega), \; \; \forall u \in \PSH(\Omega,\omega). $$

We thus define a generalization of the usual Cegrell classes as follows.
\begin{definition}\label{Def Gen Ceg class}
 Let $\mathcal{K}(\Omega) \in \{\mathcal{F}(\Omega)$,$\mathcal{E}_p(\Omega)$, $\mathcal{N}(\Omega)$,$\mathcal{E}(\Omega)\}$. For every maximal function $H \in \mathcal{E}(\Omega)$, we define the set
$\mathcal{K}(\Omega,\omega, H)$ by 
$$ u \in  \mathcal{K}(\Omega,\omega, H) \Leftrightarrow u \in \PSH(\Omega,\omega) \; \text{and} \;  u + \rho \in \mathcal{K}(\Omega,H). $$
\end{definition}
Recall that  the set $\mathcal{K}(\Omega, H)$ is defined by 
$$ u \in \mathcal{K}(\Omega, H) \Leftrightarrow u \in \PSH(\Omega) \; \text{and} \; H \geq u \geq H + \tilde{u}, $$
for $\tilde{u} \in \mathcal{K}(\Omega)$ \cite{Ceg08,ACCH09,ACLR25}.  In comparison with \cite[Page 16]{Sal25}, we no longer assume $H \in \mathcal{C}^0(\bar\Omega)$. It is worthwhile mentioning that 
 $\mathcal{K}(\Omega,\omega, H)$ is independent of the choice of the function $\rho$ (see \cite[Remark 4.3]{Sal25}). For $H=0$, we use the notation $\mathcal{K}(\Omega,\omega):= \mathcal{K}(\Omega,\omega,H)$.

Note that, by the work of Cegrell \cite{Ceg98,Ceg04,Ceg08}, a consequence of the last definition is
$$ v \in \mathcal{K}(\Omega,\omega,H) \; \text{and} \; u \in \PSH(\Omega,\omega) \; \text{with} \; v \leq u \leq H \Longrightarrow  u \in \mathcal{K}(\Omega,\omega,H).  $$

By \cite[Proposition 4.5]{Sal25}, we know that $\mathcal{K}(\Omega,\omega,H) \subset \mathcal{E}(\Omega,\omega)$. The following result proves that the operator $(\omega + dd^c.)^n$ is well defined on the classes $\mathcal{K}(\Omega,\omega,H)$ \cite[Definition 4.9 and Theorem 4.10]{Sal25}. 
\begin{theorem}\label{thm : continuity of MA measures}
    The Monge-Ampère operator $(\omega + dd^c .)^n$ is well defined on $\mathcal{E}(\Omega,\omega)$. Furthermore, if $u_j$, $u \in \mathcal{E}(\Omega,\omega)$ are such that $u_j \searrow u$ then 
    $$ (\omega + dd^c u_j)^n \rightarrow (\omega + dd^c u)^n  $$
    in the weak sense of Radon measures. 
\end{theorem} 

Fix $u \in \PSH(\Omega,\omega)$. Fix $t$, $s\in \mathbb{R}_+$ and set $u_t = \max(u,-t)$. By \cite[Theorem 5.3]{Sal25}, we know that 
$$ \textit{1}_{\{u_{t}>-s\}}(\omega + dd^c u_{t})^n = \textit{1}_{\{u_{t}>-s\}}(\omega + dd^c \max(u_{t},-s))^n.  $$
 If we choose $t\geq s$ then 
\begin{align*}
   \textit{1}_{\{u>-t\}}(\omega + dd^c u_{t})^n &\geq 
   \textit{1}_{\{u>-s\}}(\omega + dd^c u_{t})^n \\
   &=  \textit{1}_{\{u_{t}>-s\}}(\omega + dd^c u_{t})^n \\
   &= \textit{1}_{\{u>-s\}}(\omega + dd^c u_{s})^n. 
\end{align*}
It follows that the sequence of general term
$$ \textit{1}_{\{u> -t\}} (\omega + dd^c \max(u,-t))^n, \; \; t \geq 0, $$
is non-decreasing. We define the non-pluripolar Monge-Ampère measure of $u$ by the limit
$$ NP(\omega + dd^c u)^n := \lim_{t\rightarrow +\infty} \textit{1}_{\{u> -t\}} (\omega + dd^c \max(u,-t))^n. $$
By \cite[Theorem 5.3]{Sal25}, if $u \in \mathcal{E}(\Omega,\omega)$ then  
$$ NP(\omega + dd^c u)^n = \textit{1}_{\{u> -\infty \}} (\omega + dd^c u)^n  $$ 
is the non-pluripolar part of $(\omega + dd^c u)^n$.

The following result is a generalization of \cite[Theorem 2.2]{ACLR25} and \cite[Corollary 5.5]{Sal25}. 
\begin{theorem}\label{thm Dem ineq}
    Let $u$, $v \in \mathcal{E}(\Omega,\omega)$. We have 
    $$ NP(\omega + dd^c \max(u,v))^n \geq  \textit{1}_{\{u \geq v\}} NP(\omega + dd^c u)^n +  \textit{1}_{\{u<v\}} NP(\omega + dd^c v)^n. $$
    In particular, if $u \leq v$ then 
    $$ \textit{1}_{\{u = v\}} NP(\omega + dd^c u)^n \leq \textit{1}_{\{u = v\}} NP(\omega + dd^c v)^n. $$
\end{theorem}
\begin{proof}
    For every $t \geq 0$, we set $u_t = \max(u,-t)$ and $v_t = \max(v,-t)$. Since $u_t$ and $v_t$ are bounded on $\Omega$, we get by  \cite[Corollary 5.5]{Sal25} that 
     $$ (\omega + dd^c \max(u_t,v_t))^n \geq  \textit{1}_{\{u_t \geq v_t\}} (\omega + dd^c u_t)^n +  \textit{1}_{\{u_t<v_t\}} (\omega + dd^c v_t)^n. $$
     Multiplying the last inequality by $\textit{1}_{\{\min(u,v)>-t\}}$, we get 
     by \cite[Theorem 5.4]{Sal25} that 
     \begin{align*}
         \textit{1}_{\{\min(u,v)>-t\}} (\omega + dd^c \max(u,v))^n &\geq   \textit{1}_{\{u\geq v>-t\}} (\omega + dd^c u)^n \\ &+  \textit{1}_{\{-t <u<v\}} (\omega + dd^c v)^n.
     \end{align*}
      The result then follows by letting $t \rightarrow +\infty$. 
\end{proof}

    Fix $u$, $v \in \PSH(\Omega,\omega)$.  We say that $v$ is less singular than $u$, and write  $u \preceq v$, if the difference $u - v$ is locally bounded from above in $\Omega$: that is, for every compact subset $K \subset \Omega$, there exists a constant $C_K \in \mathbb{R}$ such that  
$$ u(z) - v(z) \leq C_K \quad \text{for all } z \in K. $$
If both $u \preceq v$ and $v \preceq u$, we say that $u$ and $v$ have the same singularity type, and denote this by $u \simeq v$.

The following proposition proves that the pluripolar part of the Monge-Ampère measure of functions in $\mathcal{E}(\Omega,\omega)$ is monotone with respect to the partial preorder  $\preceq$. 
\begin{prop}\label{prop 2.5}
If $u$, $v \in \mathcal{E}(\Omega,\omega)$ are such that $u \preceq v$, then
$$ \textit{1}_{\{u=-\infty\}}(\omega + dd^c u)^n \geq \textit{1}_{\{v=-\infty\}}(\omega + dd^c v)^n. $$
\end{prop}
The last result follows directly from \cite[Lemma 4.1]{ACCH09} using Lemma \ref{lem 2.5} below. 
\begin{lemma}\label{lem 2.5}
    For every $u \in \mathcal{E}(\Omega,\omega)$, we have 
    $$ (\omega + dd^c u)^n = (dd^c (u+\rho))^n \; \; \text{on} \; \{u=-\infty\}. $$
\end{lemma}
\begin{proof}
    On one hand, by \cite[Definition 4.9]{Sal25}, we know that 
    \begin{equation}\label{eq: def MA op}
        (\omega + dd^c u)^n = \sum_{k=0}^n \binom{n}{k} (-1)^{n-k} (dd^c (u+\rho))^k \wedge (dd^c \rho - \omega)^{n-k}.
    \end{equation}
  On the other hand, the proof of \cite[Proposition 5.2]{Sal25} shows that 
    $$ (dd^c (u+\rho))^k \wedge (dd^c \rho - \omega)^{n-k}(\{u=-\infty\}) = 0, \; \; \forall k=0, ..., {n-1}. $$
    Therefore, multiplying \eqref{eq: def MA op} by $\textit{1}_{\{u=-\infty\}}$ yields the result.
\end{proof}

In the sequel, we will need the following proposition.
\begin{prop}\label{prop 2.1}
 Let $\Omega' \Subset \Omega$ be a relatively compact domain in $\Omega$. Fix $u \in \PSH(\Omega,\omega)\cap L^\infty$ and $v \in \PSH(\Omega',\omega) \cap L^\infty$. If $u \geq v$ on $\partial \Omega'$ then the function
 \[
 \varphi = 
\begin{cases}
\max(u,v) & \text{on } \; \Omega' \\
u & \text{on} \; \Omega \setminus \Omega', 
\end{cases}
\]
belongs to $\PSH(\Omega,\omega) \cap L^\infty$. 
Furthermore, if $\mu$ is a positive Radon measure on $\Omega$ such that 
$$(\omega + dd^c u)^n \geq \mu \; \text{on} \; \Omega \; \; \text{and} \; (\omega + dd^c v)^n \geq \mu \; \text{on} \; \Omega', $$
then $(\omega + dd^c \varphi)^n \geq \mu$.
\end{prop}
\begin{proof}
It is easy to see that $\varphi$ is upper semi-continuous on $\Omega$. 
    For $j \geq 1$, we set $v_j = v-1/j$ and 
    \[
 \varphi_j = 
\begin{cases}
\max(u,v_j) & \text{on } \; \Omega', \\
u & \text{on} \; \Omega \setminus \Omega'. 
\end{cases}
\]
Since $u \geq v$ on $\partial \Omega'$, we obtain $\max(u,v_j) = u$ in a neighborhood of $\partial \Omega'$, which implies $\varphi_j \in \PSH(\Omega,\omega)$. Using that $\varphi$ is upper semi-continuous and that $\varphi_j \nearrow \varphi$, we obtain $\varphi \in \PSH(\Omega,\omega)$.

Consider $\mu$ as in the proposition. According to Theorem \ref{thm Dem ineq}, we have
$$ (\omega + dd^c \varphi_j)^n \geq \mu \; \text{on} \; \Omega, \; \; \forall j \geq 1. $$
Letting $j \rightarrow +\infty$, we get by \cite[Proposition 2.1]{KN15Phong} that 
$$ (\omega + dd^c \varphi)^n \geq \mu \; \text{on} \; \Omega,$$
which finishes the proof. 
\end{proof}
\subsection{Quasi-psh envelopes}
It is well known that the method of upper psh envelopes due to Perron played a significant role in establishing the first solutions of the complex Monge-Ampère equation \cite[Section 8]{BT76}. Since then, quasi-psh envelopes have received considerable attention from many authors (see, for example, \cite{BT82,ACCH09,GLZ19,ACLR25,Sal25}, and the references therein for more information). Studying singular Monge-Ampère equations of the form \eqref{Main eq} leads us  to consider generalized quasi-psh envelopes:
\begin{theorem}\label{thm: omega-psh env}
   Let $v \in \mathcal{E}_0(\Omega) \cap \mathcal{C}^0(\bar\Omega)$, and let $\psi \in \PSH(\Omega,\omega)\cap L^\infty$ with $\psi = 0$ on $\partial \Omega$. The function
   $$ u := \sup \{h \in \PSH(\Omega,\omega) : (\omega+ dd^c h)^n \geq (\omega+ dd^c \psi)^n, \; \text{and} \; h \leq v - \rho\} $$
   belongs to $\PSH(\Omega,\omega)\cap L^\infty$, satisfies $v+\psi \leq u \leq v - \rho$ on $\Omega$, and we have
   $$ (\omega + dd^c u)^n \geq (\omega + dd^c \psi)^n \; \; \text{on} \; \Omega $$
   with equality on $\{u < v-\rho\}$. 
\end{theorem}
The proof follows the lines of the one of \cite[Corollary 9.2]{BT82} that treats the case $\omega = 0$ and $(\omega + dd^c \psi)^n = 0$. Thus, we are lead to prove the following generalization of \cite[Proposition 9.1]{BT82}.
\begin{lemma}\label{lem: constructing upper functions}
    Fix $\psi \in \PSH(\Omega,\omega)\cap L^\infty$ with $\psi = 0$ in $\partial \Omega$. Assume $u \in \PSH(\Omega,\omega) \cap L^\infty$ is such that $(\omega + dd^c u)^n \geq (\omega + dd^c \psi)^n$. Then, for every ball $D \Subset \Omega$, there is $u_D \in \PSH(\Omega,\omega) \cap L^\infty(\Omega)$ such that
    \begin{center}
    $u_D \geq u$ on $\Omega$ with equality on $\Omega \setminus D$ 
       \end{center}       
    and 
    \begin{center}
            $(\omega + dd^c u_D)^n \geq (\omega + dd^c \psi)^n$ with equality on $D$. 
       \end{center}
\end{lemma}
\begin{proof}
The proof is divided into two steps.

\textbf{Step 1.} We first assume that 
$$ (\omega + dd^c \psi)^n \leq (dd^c  h)^n \; \; \text{on} \; D, $$
for some $h \in \mathcal{E}_0(D)$. 
Let $u_j \in \mathcal{C}^0(\partial D)$ be a decreasing sequence that converges to $u$. By \cite[Theorem 3.1]{KN23a}, we can find $v_j \in \PSH(D,\omega) \cap L^\infty$ such that $v_j = u_j$ on $\partial D$ and $(\omega + dd^c v_j)^n = (\omega + dd^c \psi)^n$ on $D$. According to \cite[Corollary 3.4]{KN15Phong}, the sequence $(v_j)$ is decreasing and $v_j \geq u$ on $D$. Setting $v = \lim v_j$, we have $v \in \PSH(D,\omega) \cap L^\infty$, $v = u$ on $\partial D$ and 
$$ (\omega + dd^c v)^n = (\omega + dd^c \psi)^n, $$
by \cite[Proposition 1.2]{KN15Phong}. 
Consider the function
\[
\varphi= 
\begin{cases}
v & \text{on } \; D, \\
u & \text{on} \; \Omega \setminus D. 
\end{cases}
\]
According to Proposition \ref{prop 2.1}, we know that $\varphi \in \PSH(\Omega,\omega)\cap L^\infty$ and that  $(\omega + dd^c \varphi)^n \geq (\omega + dd^c \psi)^n$, hence we can take $u_D = \varphi$.

    \textbf{Step 2.} We now move on to the general case. Fix $h \in \mathcal{E}_{0}(D)$ and set 
    $$ \mu_j = \textit{1}_{D \cap \{\psi + \rho > j h\} \cup (\Omega \setminus D)} (\omega + dd^c \psi)^n. $$
    Observe that $\mu_j$ is a positive Radon measure on $\Omega$ and that $\mu_j \nearrow (\omega + dd^c \psi)^n$ on $\Omega$ as $j \rightarrow +\infty$. Moreover,  Theorem \ref{thm Dem ineq} implies
    $$ \mu_j \leq (dd^c \max(\rho + \psi, jh))^n \; \; \textit{on} \; D. $$
    Since $\max(\rho + \psi, jh) \in \mathcal{E}_0(D)$, Step 1 gives 
  $u_D^j \in \PSH(\Omega,\omega) \cap L^\infty$ which
    satisfies the conclusion of the lemma with $(\omega + dd^c \psi)^n$ is replaced by $\mu_j$. By \cite[Proposition 1.2 and Corollary 3.4]{KN15Phong}, the sequence $(u_D^j)$ decreases on $\Omega$ to a function $u_D$ which has all the required properties.
\end{proof}
We now proceed to the proof of Theorem \ref{thm: omega-psh env}.
\begin{proof}[Proof of Theorem \ref{thm: omega-psh env}]
We first prove $u \in \PSH(\Omega,\omega) \cap L^\infty$. Observe that $u \geq v+\psi$. By Choquet's lemma \cite[Lemma 9.3]{GZ17}, there is $u_j \in  \PSH(\Omega,\omega) \cap L^\infty$ such that $u_j \leq v-\rho$, $(\omega + dd^c u_j)^n \geq  (\omega + dd^c \psi)^n$ and $(\sup_j u_j)^* = u^*$\footnote{ Here $u^*$ refers to the upper semi-continuous regularization of $u$: the smallest upper semi-continuous function above $u$.}. By Theorem \ref{thm Dem ineq}, up to replace $u_j$ by $\max(u_1,...,u_j)$, we can suppose $(u_j)$ is non-decreasing. That implies $u_j$ converges a.e. to $u^*$, thus $u^* \leq v - \rho$ everywhere on $\Omega$ since these are quasi-psh functions. Also, \cite[Proposition 2.1]{KN15Phong} implies $(\omega + dd^c u^*)^n \geq (\omega + dd^c \psi)^n$, hence $u = u^* \in \PSH(\Omega,\omega) \cap L^\infty$.

It then remains to prove that $(\omega + dd^c u)^n = (\omega + dd^c \psi)^n$ on $\{u < v-\rho\}$. Fix $z \in \{u<v-\rho\}$. By continuity, we can find a small ball $D$ centered on $z$ such that  $\sup_{\bar D} (u+\rho) < \inf_{\bar D} v$. It suffices to prove $(\omega + dd^c u)^n = (\omega + dd^c \psi)^n$ on $D$. By the last lemma, there is $u_D \in \PSH(\Omega,\omega) \cap L^\infty$ such that $u_D \geq u$ on $\Omega$, $u_D = u$ on $\Omega \setminus D$ and $(\omega + dd^c u_D)^n \geq (\omega + dd^c \psi)^n$ on $\Omega$ with equality on $D$. On one hand, we have by the above $u_D \leq v -\rho$ on $\Omega \setminus D$. On the other hand, we have 
$$ \sup_{\bar D} (u_D + \rho) = \sup_{\partial D} (u + \rho)< \inf_{\bar D} v, $$
hence $u_D \leq v - \rho$ on $\Omega$, thus $u_D = u$. The proof is complete. 
\end{proof}
\section{Solving degenerate Monge-Ampère equations}\label{sec: solving MAE}
Fix $H \in \mathcal{E}(\Omega)\cap \MPSH(\Omega)$ and let $\mathcal{K}(\Omega,\omega,H)$  denote a generalized Cegrell class. 
The objective of this section is to study the existence of weak solutions to the equation
\begin{equation}\label{MA equation}
    (\omega + dd^c \varphi)^n = \mu, \; \; \varphi \in \mathcal{K}(\Omega,\omega,H),
\end{equation}
where $\mu$ is a positive Radon measure $\Omega$. 

In a previous work, by assuming $H \in \mathcal{C}^0(\bar\Omega)$ and $\mu$ vanishes on pluripolar sets,  the last author proved the existence of solutions to \eqref{MA equation} provided that $\mu$ is dominated by the Monge-Ampère measure of a function in $\mathcal{K}(\Omega)$. 

Here we prove the following.
\begin{theorem}\label{main thm}
    Assume $\mu \leq (dd^c u)^n$ is dominated by the Monge-Ampère measure of a function $u \in \mathcal{K}(\Omega)$. Then, there is $\varphi \in \mathcal{K}(\Omega,\omega,H)$ such that 
    $$ (\omega + dd^c \varphi)^n = \mu. $$
\end{theorem}
The proof of the last theorem is based on the following lemma.
\begin{lemma}\label{main lem}
    Let $\psi \in \PSH(\Omega,\omega) \cap L^\infty$ with $\psi = 0$ on $\partial \Omega$. Let $v \in \mathcal{E}(\Omega)$ be such that $(dd^c v)^n$ is carried by a pluripolar set. Then, the function 
    $$ \varphi = \sup \{h \in \mathcal{E}(\Omega,\omega) : (\omega + dd^c  h)^n \geq (\omega + dd^c \psi)^n, \; \text{and} \; h \leq v -\rho\} $$
    belongs to $\mathcal{E}(\Omega,\omega)$, satisfies $\psi + v \leq \varphi \leq v - \rho$, and we have 
    $$ (\omega + dd^c \varphi)^n = (\omega + dd^c \psi)^n + (dd^cv)^n. $$
\end{lemma}
\begin{proof}
By \cite[Theorem 2.1]{Ceg04}, there is $v_j \in \mathcal{E}_0(\Omega) \cap \mathcal{C}^0(\bar\Omega)$ such that $v_j \searrow v$. Consider
$$ \varphi_j = \sup \{ h \in \mathcal{E}(\Omega,\omega) : (\omega + dd^c  h)^n \geq (\omega + dd^c \psi)^n, \; \text{and} \; h \leq v_j -  \rho \}. $$ 
By Theorem \ref{thm: omega-psh env}, we know that $\varphi_j \in \PSH(\Omega,\omega) \cap L^\infty$ and $v_j + \psi \leq \varphi_j \leq v_j - \rho$.
Since $\varphi_j \searrow \varphi$, we get $\varphi \in \mathcal{E}(\Omega,\omega)$ and $v+\psi \leq \varphi \leq v-\rho$. Again, by Theorem \ref{thm: omega-psh env}, we have
 $$ (\omega + dd^c \varphi_j)^n =  (\omega + dd^c \psi)^n + \textit{1}_{\{\varphi_j=v_j-\rho\}} (\omega + dd^c \varphi_j)^n. $$ 
On one hand, Theorem \ref{thm Dem ineq} gives 
$$ (\omega + dd^c {\varphi_j})^n \leq (\omega + dd^c \psi)^n + \textit{1}_{\{\varphi_j = v_j - \rho\}} (dd^c (\rho+\varphi_j))^n \leq (\omega + dd^c \psi)^n + (dd^c v_j)^n. $$
Taking $j\rightarrow +\infty$, we then obtain by Theorem \ref{thm : continuity of MA measures} that
\begin{equation}\label{MA ineq}
(\omega + dd^c \psi)^n \leq (\omega + dd^c \varphi)^n \leq (\omega + dd^c \psi)^n + (dd^c v)^n.    
\end{equation}
On the other hand, we have $v + \psi \leq \varphi \leq v -\rho$, hence  $\varphi - v$ is bounded in $\Omega$. It follows from Proposition \ref{prop 2.5} and Lemma \ref{lem 2.5} that 
$$ \textit{1}_{\{\varphi=-\infty\}} (\omega + dd^c \varphi)^n = \textit{1}_{\{\varphi=-\infty\}}(dd^c (\varphi+\rho))^n = \textit{1}_{\{v=-\infty\}} (dd^c v)^n = (dd^c v)^n. $$
From this and \eqref{MA ineq}, we obtain
\begin{align*}
    (\omega + dd^c \varphi)^n &= \textit{1}_{\{\varphi>-\infty\}} (\omega + dd^c \varphi)^n + \textit{1}_{\{\varphi = -\infty\}} (\omega + dd^c \varphi)^n \\
    &= \textit{1}_{\{\varphi>-\infty\}} (\omega + dd^c \psi)^n + (dd^c v)^n \\
    &= (\omega + dd^c \psi)^n + (dd^c v)^n, 
\end{align*}
finishing the proof. 
\end{proof}
We now proceed to the proof of Theorem \ref{main thm}. 
\begin{proof}[Proof of Theorem \ref{main thm}]
According to \cite[Theorem 5.11]{Ceg04}, we can write 
$$\mu = f(dd^c \psi)^n +\nu, $$
where $\psi \in \mathcal{E}_0(\Omega)$, $0 \leq f \in L^1((dd^c \psi)^n)$ and $\nu$ is a positive Radon measure carried by a pluripolar set. Applying \cite[Theorem 4.13]{ACCH09},  we can further write $\nu = (dd^c w)^n$ with 
$w \in \mathcal{K}(\Omega,H)$. By \cite[Theorem 3.1]{KN23a}, there is $\psi_j \in \PSH(\Omega,\omega) \cap L^\infty$ such that $\psi_j = 0$ on $\partial \Omega$ and 
 $$  (\omega + dd^c \psi_j)^n = \min(f,j) (dd^c \psi)^n. $$
 Moreover, \cite[Corollary 3.4]{KN15Phong} shows that $\psi_j$ is decreasing; set $\tilde{\psi} = \lim \psi_j$. The proof of \cite[Theorem 5.8]{Sal25} shows that $\tilde{\psi} \in \mathcal{K}(\Omega,\omega)$ and that
 $$  (\omega + dd^c \tilde{\psi})^n = f (dd^c \psi)^n. $$
 Consider 
 $$ \varphi_j := \sup \{h \in \mathcal{E}(\Omega,\omega) : (\omega + dd^ch)^n \geq (\omega+ dd^c \psi_j)^n \; \text{and} \; h \leq w - \rho \}.  $$
 The last lemma, applied with  $\psi$ replaced by $\psi_j$, yields $\varphi_j \in \mathcal{E}(\Omega,\omega)$, $\psi_j + w \leq \varphi_j \leq w - \rho$ and 
 \begin{equation}\label{eq MA meas varphi j}
     (\omega + dd^c \varphi_j)^n = (\omega+ dd^c \psi_j)^n + (dd^c w)^n.
 \end{equation}
 Moreover, since
 $$ (\omega+ dd^c \psi_j)^n = \min(f,j) (dd^c\psi)^n, $$
 is increasing with respect to $j$, we infer that the sequence $(\varphi_j)$ is decreasing. 
Setting $\varphi = \lim \varphi_j$, we have $\varphi \in \PSH(\Omega,\omega)$ and  $\tilde{\psi} + w \leq \varphi \leq w- \rho$, hence $\varphi \in \mathcal{K}(\Omega,\omega,H)$.  Taking the limit $j \rightarrow +\infty$ in \eqref{eq MA meas varphi j}, we get by Theorem \ref{thm : continuity of MA measures} that 
 $$ (\omega + dd^c \varphi)^n =  \mu, $$
 proving the theorem. 
\end{proof}
As a consequence of Theorem \ref{main thm}, we obtain the following characterization of the Range of the operator $(\omega + dd^c.)^n$, generalizing \cite[Corollary 5.11]{Sal25}.
\begin{corollary}
  Let $\mu$ be positive Radon measure on $\Omega$. The following affirmations are equivalent. 
  \begin{itemize}
      \item[(i)] The equation $(dd^c u)^n = \mu$ admits a solution $u \in \mathcal{K}(\Omega)$.
       \item[(ii)]  The equation $(\omega + dd^c \varphi)^n = \mu$ has a solution $\varphi \in \mathcal{K}(\Omega,\omega)$.
       \item[(iii)]  For every $H \in \mathcal{E}(\Omega) \cap \MPSH(\Omega)$, there is $\varphi \in \mathcal{K}(\Omega,\omega,H)$ such that $(\omega + dd^c \varphi)^n = \mu$.
  \end{itemize}
\end{corollary}
\begin{proof}
    The implication (i) $\Rightarrow$ (iii) follows from Theorem \ref{main thm} and (ii) is a particular case of (iii).  An application of \cite[Theorem 4.13]{ACCH09} gives (ii) $\Rightarrow$ (i), which finishes the proof. 
\end{proof}
The following result studies the uniqueness of solutions to \eqref{MA equation}.  
\begin{theorem}\label{thm : uniqueness}
Let $H \in \mathcal{E}(\Omega)\cap \MPSH(\Omega)$, and let $u$, $v \in \mathcal{N}(\Omega,\omega,H)$ be such that $u \preceq v$. 
\begin{itemize}
    \item[(1)] If $(\omega + dd^c u)^n \leq ( \omega + dd^c v)^n$  then $u \geq v$.
    \item[(2)]  If $(\omega + dd^c u)^n = ( \omega + dd^c v)^n$  then $u = v$.
\end{itemize}
\end{theorem}
\begin{proof}
We start by proving (1). 
Let $\tilde{u} \in \mathcal{N}(\Omega)$ be such that $\tilde{u} +H \leq u + \rho \leq H$.  Observe that $u - v = u + \rho - (v+\rho) \geq \tilde{u}$.
Therefore, setting 
$$ \varphi :=  \sup \{ h \in \PSH(\Omega)\; ; \; h \leq \min(u - v,0) \},  $$
we have $\varphi \in \PSH(\Omega)$ and $\tilde{u} \leq \varphi \leq 0$, hence $\varphi \in \mathcal{N}(\Omega)$ by \cite[Corollary 3.14]{Ceg08}.  Moreover, Proposition \ref{prop 2.5} together with the hypothesis imply  
$$ \textit{1}_{\{u=-\infty\}} (\omega + dd^c u)^n = \textit{1}_{\{v=-\infty\}} (\omega + dd^c v)^n.  $$
Hence 
$$ \textit{1}_{\{u=-\infty\}} (dd^c(u+\rho))^n = \textit{1}_{\{v=-\infty\}} (dd^c(v+\rho))^n  $$
by Lemma \ref{lem 2.5}. It thus follows from \cite[Lemma 3.3]{ACLR25} that $(dd^c\varphi)^n$ does not charge pluripolar sets; in particular $NP(dd^c \varphi)^n = (dd^c \varphi)^n$.

Consider now $D = \{\varphi = \min(u - v,0)\}$. We know by \cite[Theorem 2.7]{ACLR25} that $(dd^c \varphi)^n$ is carried by $D$. Since $\varphi + v \leq u$, Theorem \ref{thm Dem ineq} implies that
\begin{align*}
    (dd^c \varphi)^n + \textit{1}_D NP(\omega + dd^c v)^n &= (dd^c \varphi)^n + \textit{1}_{D \cap \{v>-\infty\}} (\omega + dd^c v)^n \\
    &= \textit{1}_{D \cap \{v>-\infty\}} \left( (dd^c \varphi)^n + (\omega + dd^c v)^n \right) \\
    &\leq \textit{1}_D NP(\omega + dd^c u)^n \\
    &\leq \textit{1}_D NP(\omega + dd^c v)^n,
\end{align*}
which proves $(dd^c \varphi)^n = 0$. It then follows from \cite[Theorem 3.6]{Ceg08} that $\varphi = 0$, and thus $u \geq v$. 

The second assertion follows directly from the first. Indeed, the affirmation (1) implies $u \geq v$, whence $u \simeq v$ by hypothesis. Therefore, changing the role of $u$ and $v$ in (1) yields $u = v$. 
\end{proof}

\pagestyle{empty}

\begin{thebibliography}{BCHM10}
\bibliographystyle{plain}

\bibitem[{\AA}CCH09]{ACCH09}
P.~{\AA}hag, U.~Cegrell, R.~Czy{\.z}, and P. H. Hiep,
  \emph{{M}onge--{A}mp{è}re measures on pluripolar sets}, Journal de
  math{é}matiques pures et appliqu{é}es \textbf{92} (2009), no.~6,
  613--627.

\bibitem[{\AA}CLR25]{ACLR25}
P.~{\AA}hag, R.~Czy{\.z}, C.H. Lu, and A. Rashkovskii,  \emph{Geodesic connectivity and rooftop envelopes in the Cegrell classes}, Math. Ann. 391, 3333–3361 (2025). https://doi.org/10.1007/s00208-024-03003-7
  

\bibitem[BT76]{BT76}
E. Bedford and B. A. Taylor, \emph{The Dirichlet problem for a complex Monge-Ampère operator,}  Invent. Math. 37, 1–44 (1976). 

\bibitem[BT82]{BT82}
E. Bedford and B. A. Taylor, \emph{A new capacity for plurisubharmonic functions,} Acta. Math. 149, 1–40 (1982). 

\bibitem[Ceg98]{Ceg98}
 U. Cegrell, \emph{Pluricomplex energy,} Acta Math, 180 (1998), 187-217.  DOI: 10.1007/BF02392899

\bibitem[Ceg04]{Ceg04}
U. Cegrell, \emph{The general definition of the complex Monge-Ampère operator,} Ann. Inst. Fourier (Grenoble) 54 (2004), 159-179. (www.numdam.org/articles/10.5802/aif.2014/)


\bibitem[Ceg08]{Ceg08}
 U. Cegrell, \emph{A general Dirichlet problem for the
complex Monge-Ampère operator,} Ann. Polon.  Math.
94.2 (2008) 131-147. (http://eudml.org/doc/280684)


\bibitem[GLZ19]{GLZ19}
 V. Guedj, C. H. Lu and A. Zeriahi, \emph{Plurisubharmonic envelopes and supersolutions,} Journal of  Differential Geometry
113 (2019) 273-313.


\bibitem[GZ17]{GZ17}
 V. Guedj and  A. Zeriahi, \emph{Degenerate complex Monge–Ampère equations,} EMS Tracts. Math. 26 (2017)

\bibitem[HHHP14]{HHHP14} 
L. M. Hai,  P. H. Hiep, N. X. Hong, and  N. V. Phu,  \emph{The Monge–Ampère type equation in the weighted pluricomplex energy class}, International Journal of Mathematics, 25(05), 1450042. (2014).



 \bibitem[Kol95]{Kol95}
S. Ko{\l}odziej, \emph{The range of the complex Monge-Ampère operator II},
Indiana Univ. Math. J. 44 (1995), 765-782.


 \bibitem[Kol05]{Kol05}
S. Ko{\l}odziej, \emph{The complex Monge-Ampère equation and pluripotential
theory}, Mem. Amer. Math. Soc. 178 (2005).
 
 \bibitem[KN15]{KN15Phong}
S. Ko{\l}odziej and N. C. Nguyen, \emph{Weak solutions to the complex Monge-Ampère equation on Hermitian manifolds,} Analysis, complex geometry, and mathematical physics.  Contemporary Mathematics, vol. 644 (American Mathematical Society, Providence, RI, 2015) 141-158. 





\bibitem[KN23a]{KN23a}
S. Ko{\l}odziej and N. C. Nguyen, \emph{The Dirichlet problem for the Monge-Ampère equation on Hermitian manifolds with boundary,} Calc. Var. Partial Differential Equations
62 (2023), no. 1, Paper No. 1.

\bibitem[KN23b]{KN23b}
S. Ko\l odziej and N. C. Nguyen, \emph{Weak Solutions to Monge–Ampère Type Equations on Compact Hermitian Manifold with Boundary,} The Journal of Geometric Analysis 33.1 (2023): 15.

\bibitem[QH19]{QH19}
V. V. Quan and L. M.  Hai, \emph{ Weak solutions to the complex Monge-Ampère equation on open subsets of $\mathbb{C}^n$}, Matematychni Studii, 51(2)
(2019)

\bibitem[Sal25]{Sal25}
M. Salouf, \emph{Degenerate complex Monge-Ampère equations with non-K\"ahler forms in bounded domains}, Indiana Univ. Math. J. 74 No. 1 (2025), 131–156.

\bibitem[Zer97]{Zer97}
A. Zeriahi, \emph{Pluricomplex Green functions and the Dirichlet problem for the complex Monge-Ampère operator},  Michigan Mathematical Journal, Michigan Math. J. 44(3), 579-596, (1997)

\end{thebibliography}
\end{document}